\documentclass[10pt]{article}
\usepackage{latexsym}
\title{On the Minimum Width of a Cutset in the Truncated Boolean Lattice}

\author{B\'{e}la Bajnok\\ Gettysburg College\\ Gettysburg, PA 17325-1486 USA}

\newtheorem{thm}{Theorem}

\newtheorem{conj}[thm]{Conjecture}

\newcommand{\bc}[2]{{{#1}\choose{#2}}}

\begin{document}

\renewcommand{\today}{March 13, 1998}

\maketitle

\begin{abstract}

For integers $0 \leq m \leq l \leq n-m$, the truncated Boolean lattice ${\cal B}_n(m,l)$ is the poset of all subsets of $[n] = \{1, 2, \ldots, n\}$ which have size at least $m$ and at most $l$. ${\cal C} \subseteq {\cal B}_n(m,l)$ is a {\em cutset} if it meets every chain of length $l-m$ in ${\cal B}_n(m,l)$, and the {\em width} of ${\cal C}$ is the size of the largest antichain in ${\cal C}$. We conjecture that for $n >> m$ the minimum width $h_n(m,l)$ of a cutset in ${\cal B}_n(m,l)$ is $\Sigma_{j \geq 0} \Delta_n(m-jc) = \Delta_n(m)+\Delta_n(m-c)+\Delta_n(m-2c)+ \dots$, where $c=l-m+1$ is the number of level sets in ${\cal B}_n(m,l)$ and $\Delta_n(k)={n \choose k}- {n \choose k-1}$. We establish our conjecture for the cases of ``short lattices'' ($l=m$, $l=m+1$, and $l=m+2$). For ``taller lattices'' ($l \geq 2m$) our conjecture gives ${n \choose m} - {n \choose m-1}$, independently of $l$. Our main result is that $h_n(m,l) \leq {n \choose m} - {n \choose m-1}$ if $l \geq 2m$.

\end{abstract}

Let ${\cal B}_n=2^{[n]}$ be the {\em Boolean lattice} of order $n$, that is the lattice of all subsets (often called {\em nodes}) of $[n] = \{1, 2, \ldots, n \}$ ordered by inclusion. For $0 \leq i \leq n$ we define the $i$-th {\em level set} $\bc{[n]}{i}$ of $2^{[n]}$ as the set of all subsets of size $i$. The {\em truncated Boolean lattice} ${\cal B}_n(m,l)$ ($0 \leq m \leq l \leq n$) is then the union of level sets $\bc{[n]}{i}$ for $m \leq i \leq l$. Because of symmetry, throughout this paper we assume that $m \leq l \leq n-m$.   

A collection of subsets $A_1 \subset A_2 \subset \cdots \subset A_c$ in $2^{[n]}$ is called a {\em chain} of {\em size} $c$ and {\em length} $c-1$; this chain is called {\em saturated} if $|A_{i+1}|=|A_i|+1$ for every $i=1,2,\dots,c-1$. All chains in this paper will be saturated. A {\em maximal chain} in ${\cal B}_n(m,l)$ is one that has size $c=l-m+1$. A collection of $w$ nodes with the property that none of them contains another is called an {\em antichain} of size $w$. The length and the width of a collection of subsets ${\cal A} \subseteq 2^{[n]}$ are defined as the length of the longest chain and the size of the largest antichain in ${\cal A}$, respectively. By Dilworth's Theorem, the width of ${\cal A}$ is the minimum number of chains in a chain decomposition of ${\cal A}$.

A {\em cutset} in ${\cal B}_n(m,l)$ is defined as a collection of subsets ${\cal C} \subseteq {\cal B}_n(m,l)$ which intersects all maximal chains. Let us denote the minimum width of a cutset in ${\cal B}_n(m,l)$ by $h_n(m,l)$. In \cite{BS-1} we proved that $h_n(1,n-1)=n-1$. This paper is an attempt to generalize this result. In particular, we state the following conjecture.

\begin{conj} \label{conject}

For integers $n$, $m$, $l$, and $k$ which satisfy $0 \leq m \leq l \leq n-m$, we set $c=l-m+1$ and $\Delta_n(k)=\bc {n}{k}- \bc{n}{k-1}$. Then for $n >> m$ the minimum width of a cutset in ${\cal B}_n(m,l)$ is $h_n(m,l)=\Sigma_{j \geq 0} \Delta_n(m-jc) = \Delta_n(m)+\Delta_n(m-c)+\Delta_n(m-2c)+ \dots$.

\end{conj} 

Related to $h_n(m,l)$ is $g_n(m,l)$, defined as the smallest value of $k$ for which there exists a cutset with at most $k$ nodes at each level in ${\cal B}_n(m,l)$. In \cite{BS-2} we established the values of $g_n(m,l)$     for the cases $l=m$, $l=m+1$, and $l=m+2$. Namely, after a general characterization of cutsets in terms of the number and sizes of their elements based on the Kruskal-Katona Theorem, we proved the following.

\begin{thm} \label{g-exact} \cite{BS-2} Let $m$ and $n$ be non-negative integers. Then

\begin{enumerate}

\item

$g_n(m,m) = \bc{n}{m}$ for $n \geq m$.

\item

$g_n(m,m+1) = \bc{n-1}{m}$ for $n \geq m+1$.

\item

$g_n(m,m+2) = \sum_{j=0}^{m} \bc{n-2j-2}{m-j}$  for $n \geq 2m+2$.

\end{enumerate}

\end{thm}

It is obvious that $h_n(m,l) \geq g_n(m,l)$. We state the following surprising conjecture.

\begin{conj} \label{g=h}

For integers $n$, $m$, and $l$ for which $0 \leq m \leq l \leq n-m-1$ and $n>>m$, we have $g_n(m,l)=h_n(m,l)$. In addition, $g_n(m,n-m)=\bc {n-1}{m}- \bc{n-1}{m-1}$.

\end{conj}

We can now prove the cases of $c=1$, $c=2$, and $c=3$ (i.e. $l=m$, $l=m+1$, and $l=m+2$) of Conjecture \ref{conject} and Conjecture \ref{g=h}.

\begin{thm} \label{exact} Let $m$ and $n$ be non-negative integers. Then

\begin{enumerate}

\item

$h_n(m,m) = \Sigma_{j \geq 0} \Delta_n(m-j) = \bc{n}{m}$ for $n \geq m$.

\item

$h_n(m,m+1) = \Sigma_{j \geq 0} \Delta_n(m-2j) = \bc{n-1}{m}$ for $n \geq m+1$.

\item

$h_n(m,m+2) = \Sigma_{j \geq 0} \Delta_n(m-3j) = \sum_{j=0}^{m} \bc{n-2j-2}{m-j}$  for $n \geq 2m+2$.

\end{enumerate}

\end{thm}

{\em Proof.} We leave it to the reader to establish the binomial identities $\Sigma_{j \geq 0} \Delta_n(m-j) = \bc{n}{m}$, $\Sigma_{j \geq 0} \Delta_n(m-2j) = \bc{n-1}{m}$, and $\Sigma_{j \geq 0} \Delta_n(m-3j) = \sum_{j=0}^{m} \bc{n-2j-2}{m-j}$. 

The case $h_n(m,m) = \bc{n}{m}$ is obvious. Since $h_n(m,l) \geq g_n(m,l)$, by Theorem \ref{g-exact} it suffices to prove $h_n(m,m+1) \leq \bc{n-1}{m}$ and $h_n(m,m+2) \leq \sum_{j=0}^{m} \bc{n-2j-2}{m-j}$. 

We can easily construct $\bc{n-1}{m}$ chains in ${\cal B}_n(m,m+1)$ which form a cutset. Color the nodes of ${\cal B}_n(m,m+1)$ by two colors: black if they do not contain the element ~1, and white if they do. There are exactly $\bc{n-1}{m}$ black nodes at level $m$ and the same number of white nodes at level $m+1$, and they clearly form chains of size ~2. Finally, these $\bc{n-1}{m}$ chains form a cutset in ${\cal B}_n(m,m+1)$, since the nodes that are disjoint from them are the white nodes at level $m$ and the black nodes at level $m+1$, but these nodes form an antichain in ${\cal B}_n(m,m+1)$.

For the case $l=m+2$ we first color the nodes of ${\cal B}_n(m,m+2)$ by four colors: black if they do not contain the elements ~1 and ~2, blue if they contain ~1 but not ~2, red if they contain ~2 but not ~1, and white if they contain both ~1 and ~2. Note that as we go up on a chain, the color of the elements ``fade.'' More precisely, black nodes can only be extended up to black, blue, or red nodes, blue or red nodes can be extended to the same color or white, and white nodes can only be extended to white nodes.

The number of black nodes at level $m$, blue nodes at level $m+1$, and white nodes at level $m+2$ are all exactly $\bc{n-2}{m}$, and these nodes clearly form the same number of chains of size ~3 in ${\cal B}_n(m,m+2)$. At this point in the construction, if a chain of size ~3 is disjoint from our chosen $\bc{n-2}{m}$ chains in ${\cal B}_n(m,m+2)$, it must consist entirely of red nodes. The set of red nodes is isomorphic to ${\cal B}_{n-2}(m-1,m+1)$, hence we see that $h_n(m,m+2) \leq \bc{n-2}{m}+h_{n-2}(m-1,m+1)$, from which our claim follows by induction. $\quad \Box$

Let us return to Conjecture \ref{conject}. The sum in Conjecture \ref{conject} has $\lfloor m/c \rfloor +1$ positive terms; in particular, only the first term is positive if $c>m$. According to this conjecture, $h_n(m,l)$ decreases at a rapidly decreasing rate as $l$ increases from $m$ to $n-m$ and, in fact, it becomes the constant $\bc{n}{m} - \bc{n}{m-1}$ for $2m \leq l \leq n-m$. Therefore, surprisingly, for $n>>m$ the minimum width of a cutset in ${\cal B}_n(m,2m)$ is the same as it is in ${\cal B}_n(m,n-m)$.

Our main result is the following theorem.

\begin{thm} \label{main} If $n$, $m$, and $l$ are integers such that $0 \leq 2m \leq l \leq n-m$, then $h_n(m,l) \leq \bc{n}{m} - \bc{n}{m-1}$.

\end{thm}

{\em Proof.} First note that if $2m \leq l$ then $h_n(m,l) \leq h_n(m,2m)$, so it suffices to construct $\bc{n}{m}-\bc{n}{m-1}$ chains in ${\cal B}_n(m,2m)$ which form a cutset.

We start by coloring the nodes of $2^{[n]}$ as follows. The color $c(A)$ of $A \in 2^{[n]}$ is defined as the set $A \cap [2m]$. Each color class contains $2^{n-2m}$ elements and is isomorphic to the lattice $2^{[n-2m]}$. Furthermore, the set of $2^{2m}$ color classes forms a lattice ${\cal Q}$ which is isomorphic to $2^{[2m]}$. 

Next, we consider the following recursive construction of Griggs (see \cite{Griggs}, 264-266) to cover the nodes of the Boolean lattice using chains of length at most $c$. This is trivial for $2^{[1]}$. Once the partition of $2^{[k]}$ into chains of length at most $c$ is given, we partition $2^{[k+1]}$ as follows. For every chain $A_1 \subset A_{2} \subset \cdots \subset A_r$ used in the partition of $2^{[k]}$, in $2^{[k+1]}$ we form the chain(s)

\begin{itemize}

\item

$A_1 \subset A_1 \cup \{k+1\}$ if $r=1$,

\item

$A_1 \subset A_{2} \subset \cdots \subset A_r \subset A_r \cup \{k+1\}$ and $A_1 \cup \{k+1\} \subset A_{2} \cup \{k+1\} \subset \cdots \subset A_{r-1} \cup \{k+1\}$ if $2 \leq r \leq c-1$, and

\item

$A_1 \subset A_{2} \subset \cdots \subset A_r$ and $A_1 \cup \{k+1\} \subset A_{2} \cup \{k+1\} \subset \cdots \subset A_r \cup \{k+1\}$ if $r=c$.
 
\end{itemize}
This recursive construction is the ``quit while we're ahead'' modification of de Bruijn's well known construction for a symmetric chain decomposition of the Boolean lattice. In particular, it provides a chain partition of ${\cal Q} \cong 2^{[2m]}$ using $\bc{2m}{m}$ chains with the following properties.

\begin{enumerate}

\item

Every chain has length at most $m$;

\item

Exactly $\bc{2m}{j}-\bc{2m}{j-1}$ of the chains start at level $j$ for $j=0,1, \dots m$; and

\item

If $A_1 \subset A_2 \subset \cdots \subset A_{r}$ and $B_1 \subset B_2 \subset \cdots \subset B_{s}$ are two of the chains and $A_i \subset B_j$ for some $i=1, \dots, r$ and $j=1, \dots, s$, then $i \leq j$.

\end{enumerate}

Now we are ready to describe the $\bc{n}{m}-\bc{n}{m-1}$ chains which form a cutset in ${\cal B}_n(m,2m)$. Suppose $C_j \subset C_{j+1} \subset \cdots \subset C_{k}$ is one of the chains used to cover ${\cal Q}$ above, and suppose that $|C_i|=i$ for $j \leq i \leq k$. From this chain we construct $\bc{n-2m}{m-j}$ chains in ${\cal B}_n(m,2m)$ as follows. 

For each $i=j,j+1, \dots, k$, take the $\bc{n-2m}{m-j}$ nodes of color $C_i$ at level $m+i-j$ in ${\cal B}_n(m,m+k-j)$. Since $C_j \subset C_{j+1} \subset \cdots \subset C_{k}$ is a chain and $m+k-j \leq 2m$ by property 1 above, these nodes clearly form $\bc{n-2m}{m-j}$ chains from level $m$ to level $m+k-j$ in ${\cal B}_n(m,2m)$. We repeat this procedure with the other chains in ${\cal Q}$, and, by property 2 above, get $$\sum_{j=0}^m \left( \bc{2m}{j}-\bc{2m}{j-1} \right) \bc{n-2m}{m-j}=\bc{n}{m}-\bc{n}{m-1}$$ chains in ${\cal B}_n(m,2m)$. 

It remains to show that the collection ${\cal A}$ of these chains forms a cutset. Let $ D_m \subset D_{m+1} \subset \cdots \subset D_{2m} $ be an arbitrary chain from level $m$ to level $2m$ in ${\cal B}_n(m,2m)$. We shall prove that there is a $j=0,1, \dots ,m$, for which $D_{m+j} \in {\cal A}$. The (not necessarily distinct) colors $c(D_m) \subseteq c(D_{m+1}) \subseteq \cdots \subseteq c(D_{2m})$ of these $m+1$ nodes are all covered by chains in ${\cal Q}$, therefore all the corresponding color classes intersect ${\cal A}$ at some level in ${\cal B}_n(m,2m)$. Suppose that for $j=0,1, \dots, m$, the color class of $c(D_{m+j})$ intersects ${\cal A}$ at level $l_j$. Since $m \leq l_0 \leq l_1 \leq \cdots \leq l_m \leq 2m$ by property 3 above, there must be a $j=0,1, \dots ,m$ for which $l_j=m+j$. For this $j$, $D_{m+j} \in {\cal A}$, as claimed. $ \quad \Box$


\begin{thebibliography}{99}

\bibitem{BS-1} B. Bajnok and S. Shahriari, Long Symmetric Chains in the Boolean Lattice, \emph{J. Combin. Theory Ser. A}, {\bf ~75} (1996), ~44-54.

\bibitem{BS-2} B. Bajnok and S. Shahriari, On Uniform $f$-vectors of Cutsets in the Truncated Boolean Lattice, to appear.

\bibitem{Griggs} J. R. Griggs, Saturated Chains of Subsets and a Random Walk, \emph{J. Combin. Theory Ser. A}, {\bf ~47} (1988), ~262-283.


\end{thebibliography}
\end{document}